%% file: Main_New_Format.tex
\documentclass[11pt]{article}
\pdfoutput=1
\usepackage{amsfonts,amsmath,amssymb,amsthm}
\usepackage[pdftex]{graphicx}
\usepackage[hmargin=1in,vmargin=1in]{geometry}
\usepackage[sort&compress]{natbib}				
\usepackage{booktabs}
\usepackage{setspace}
\usepackage{appendix}
\usepackage{multirow}
\usepackage[colorlinks,citecolor=blue]{hyperref}
\usepackage{enumerate}
\usepackage{caption}
\usepackage[finnish, british]{babel}
\usepackage{capt-of}
\usepackage{hyperref}
\usepackage{float}
\usepackage[flushleft]{threeparttable}
\def\qed{\rule{2mm}{2mm}}

\let\footnote=\endnote

\usepackage{tikz,graphicx,pgfplots,pgfkeys,subcaption}   

\def\addlegendimage{\csname pgfplots@addlegendimage\endcsname}

\usepackage[pagewise,mathlines]{lineno}
\mathchardef\dash="2D

\newtheorem{theorem}{Theorem}[section]

\newtheorem{condition}{Condition}[section]

\newtheorem{proposition}{Proposition}[section]
\newtheorem{assumption}{Assumption}[section]
\theoremstyle{definition}
\newtheorem{mydef}{Defininition}[section]

\newtheorem{remark}{Remark}[section]
\usepackage{etoolbox} 
\AtEndEnvironment{remark}{~\qed}%
\AtEndEnvironment{example}{~\qed}%

\begin{document}

\author{
Max Tabord-Meehan\footnote{I am grateful for Ivan Canay's advice and encouragement. I would also like to thank two anonymous referees, Peter Aronow, Eric Auerbach, Joel Horowitz, Vishal Kamat, Lilla Orr, Susan Ou, and seminar participants at Northwestern University for helpful comments. This research was supported in part through the computational resources and staff contributions provided for the Social Sciences Computing Cluster (SSCC) at Northwestern University. All mistakes are my own.
}\\
Departments of Economics\\
Northwestern University\\
\url{mtabordmeehan@u.northwestern.edu}
}

\title{Inference with Dyadic Data: Asymptotic Behavior of the Dyadic-Robust $t$-Statistic}

\maketitle

\begin{abstract}
This paper is concerned with inference in the linear model with dyadic data. Dyadic data is data that is indexed by pairs of ``units'', for example trade data between pairs of countries. Because of the potential for observations with a unit in common to be correlated, standard inference procedures may not perform as expected. We establish a range of conditions under which a $t$-statistic with the dyadic-robust variance estimator of \cite{fafgub2007} is asymptotically normal. Using our theoretical results as a guide, we perform a simulation exercise to study the validity of the normal approximation, as well as the performance of a novel finite-sample correction. We conclude with guidelines for applied researchers wishing to use the dyadic-robust estimator for inference.\end{abstract}

\noindent KEYWORDS: Regression, Dyadic Data, Dependence, t-Test, Robust Variance Estimators, Degrees of Freedom Correction.

\noindent JEL classification codes: C12, C21.

\clearpage
\section{Introduction} \label{sec:intro}
\input{introduction.tex}

\section{Setup of the Model and Asymptotic Frameworks} \label{sec:setup}
\input{setup.tex}

\section{Asymptotic Properties of $T_k$} \label{sec:results}
\input{results.tex}

\section{Simulation Evidence and a Degrees of Freedom Correction}\label{sec:simulations}
\input{simulations.tex}

\section{Conclusion}\label{sec:conclusion}
\input{conclusion.tex}

\appendix
\section{Appendix}\label{sec:appendix}
\input{appendix.tex}

\clearpage	
\bibliography{references.bib}
\nocite{*}
\clearpage

\end{document}

%% file: introduction.tex
Over the last 25 years applied microeconomics has increasingly embraced the fact that dependence in cross-sectional data can affect inference. It has been well understood since at least the work of \cite{moulton1986} that failing to account for dependence in cross-sectional studies can have dire effects. In the past, researchers explicitly modeled such dependencies and used techniques such as GLS to estimate and do inference in their models. However, modern researchers are typically not satisfied with making such strong assumptions on the dependence present in the data. It is now standard practice to account for dependence by pairing standard test statistics with so-called ``robust'' variance estimators, analogous to the heteroskedasticity-robust variance estimator of \cite{white1980}. 

In this paper we focus on inference for the regression parameters in a linear model with dyadic data. Dyadic data relates to pairs of objects; examples include data on trade between pairs of countries and data on links in a social-network setting. We will call such pairs ``dyads'' and the objects within them ``units''. Because of the paired nature of the data, dyads that share a unit in common could be correlated. In order to account for this potential dependence when conducting inference, we study the asymptotic properties of a $t$-statistic formed using a ``robust" variance estimator known in the literature as the dyadic-robust variance estimator.

 \cite{fafgub2007} were the first to propose the dyadic-robust variance estimator, under the following assumption: dyads that do \emph{not} share a unit are uncorrelated, but otherwise the dependence between dyads is unspecified. To draw an analogy with cluster-robust inference \citep[see][for an extensive survey]{cameronmiller2015}, the dyadic-dependence assumption results in an ``overlapping-cluster'' configuration of the data, with each unit defining its own cluster. Subsequently, many applied papers in economics and political science have employed the dyadic-robust estimator under this same assumption \citep[an incomplete list includes][]{aker2010, leblang2010, baldwin2012, comola2014, egger2015, lustig2015, echevarria2016, poznansky2016}. Many empirical papers with dyadic data also make reference to the dependence assumption we describe, but then compute two-way clustered standard errors as described in \cite{cammil2011}. However, \cite{cammil2014} point out that this does not account for all the potential dependencies in the dyadic setting.

We present formal results under which a $t$-statistic that uses the dyadic-robust variance estimator is asymptotically normal. Using a central limit theorem for dependency graphs proved in \cite{janson1988} and careful bounding arguments, we establish a range of assumptions under which asymptotic normality holds. We then use our results to guide a simulation study of the accuracy of a normal approximation in finite samples, and discover an important setting where such an approximation is inadequate. With these insights, we propose a novel degrees of freedom correction to help alleviate the issue, and assess the performance of this correction in simulations. 

 \cite{fafgub2007} motivate the dyadic-robust variance estimator as an extension of the spatial HAC estimator of  \cite{conley1999}. Despite Conley's work being the initial motivation, neither consistency of the dyadic-robust variance estimator nor asymptotic normality of the resulting $t$-statistic, under their maintained assumptions, follow from his results. Recently, \cite{cammil2014} have proposed the use of the dyadic-robust variance estimator in the analysis of trade data, and present simulation evidence to assess its performance.  Both \cite{fafgub2007} and \cite{cammil2014} implicitly assume an asymptotic normality result for the dyadic-robust $t$-statistic in their analysis, but do not provide conditions under which such a result may hold.  \cite{aronsam2015} prove the consistency of the dyadic-robust variance estimator for cross-sectional and panel data under more strict assumptions than those considered here, but do not attempt to study the use of this estimator for inference: although they derive the asymptotic variance of the $t$-statistic, they do \emph{not} characterize its asymptotic distribution, specifically, they do not establish conditions under which the $t$-statistic is asymptotically normal. Our paper is the first to provide a theoretical grounding for the use of a normal approximation to the dyadic-robust $t$-statistic for inference in the linear model. 

The remainder of the paper is organized as follows: In Section 2, we set up the model and the asymptotic frameworks we will study. Section 3 presents our results about asymptotic normality of the $t$-statistic. In Section 4 we study the finite-sample behavior of our approximation in a simulation study, and propose a degrees of freedom correction guided by our results. Section 5 concludes.

%% file: setup.tex
\subsection{The Model}
We will now formally describe the model. Consider a collection of $G$ units indexed by $g = 1,...,G$. The data we consider is indexed by pairs of units $(g,h)$, which we call dyads. We do not require that each possible pair of units form a dyad. Pairs of units $(g,h)$ for which a dyad exists map into dyadic indices $n = 1,2,...,N$ through the index function $n(g,h)$, where for simplicity we make the assumption that $n(g,h) = n(h,g)$ (i.e. that we treat the dyads as non-directional), and the assumption that there are no elements of the type $n(g,g)$ (i.e. that we consider only pairs between distinct units). Given a dyadic index $n$, we also define the inverse correspondence $\psi$ so that $\psi(n(g,h)) = \{g,h\}$. The model we consider is the linear model:
\begin{equation}\label{eq:linmod}
y_{n(g,h)} = \beta '{\bf x}_{n(g,h)} + u_{n(g,h)}~,
\end{equation}
where ${\bf x}_n$ is $K$-dimensional, with the standard conditions that $E[{\bf x}_nu_n] = {\bf 0}$ and $E[{\bf x}_n{\bf x}_n'] > 0$. Our focus is on performing inference on the components of the regression parameter $\beta$.

Next, we present the dependence structure we will consider. Intuitively, we want observations that do not share a unit in common to be independent, but to allow correlation between observations otherwise.  The typical assumption stated in the literature \citep[see for example][]{cammil2014, aronsam2015} is that
$$E[u_nu_{m}|{\bf x}_n,{\bf x}_{m}] = 0, \hspace{2mm} \text{unless} \hspace{2mm} \psi(n) \cap \psi(m) \ne \emptyset~.$$
Although this assumption somewhat captures the intuition presented above, we will need to sharpen it considerably in order to prove formal results about our model. We impose the following dependence assumption on the data: 

\begin{assumption}\label{ass:dep}
$\{({\bf x}_n,u_n)\}_{n=1}^N$ are identically distributed. For any two disjoint subsets $S_1$, $S_2$ of $\{1,2,...,N\}$, $\{({\bf x}_n,u_n)\}_{n\in S_1}$ is independent of $\{({\bf x}_{m},u_{m})\}_{m\in S_2}$ if $\psi(n) \cap \psi(m) = \emptyset$ for every pair $n,m$ of dyads such that $n \in S_1 ,m \in S_2$. 
\end{assumption}

\begin{remark}\label{rem:panels}
Although we focus on the setting where our data is cross-sectional and the dyads are non-directional, our analysis covers settings with directional dyads, as well as panel-data with a \emph{finite} number of time periods. For example, consider the following panel-data version of our model:
$$y_{(g,h)t} = \beta'{\bf x}_{(g,h)t} + \gamma_{g} + \gamma_{h} + \alpha_{gh} + u_{(g,h)t}~,$$
where we have explicitly indexed observations by their units, and now observations are indexed by pairs of units $(g,h)$ as well as by time $t = 1,...,T$. Note that we include $\gamma_g$ and $\gamma_h$, which are unit-level fixed effects, as well as a dyad-level fixed effect $\alpha_{gh}$. 
Let $\ddot{y}_{(g,h)t}$, ${\bf \ddot{x}}_{(g,h)t}$, and $\ddot{u}_{(g,h)t}$ denote the random variables that result from performing a \emph{within} transformation:
$$\ddot{y}_{(g,h)t} := y_{(g,h)t} - \frac{1}{T}\sum_{s=1}^Ty_{(g,h)s}~,$$
and similarly for ${\bf \ddot{x}}_{(g,h)t}$ and $\ddot{u}_{(g,h)t}$. Then the transformed model
$$\ddot{y}_{(g,h)} = {\bf \ddot{x}}_{(g,h)}\beta + \ddot{u}_{(g,h)}~,$$
where $\ddot{y}_{(g,h)}$ and $\ddot{u}_{(g,h)}$ are $T \times 1$ stacked vectors and ${\bf \ddot{x}}_{(g,h)}$ is a $T \times K$ stacked matrix, can be studied completely analogously to Model (\ref{eq:linmod}) above, with the assumptions $E[{\bf \ddot{x}}_{(g,h)}'\ddot{u}_{(g,h)}] = {\bf 0}$ and $E[{\bf \ddot{x}}_{(g,h)}'{\bf \ddot{x}}_{(g,h)}] > 0$ \citep[these assumptions are implied by standard primitive conditions on the original, untransformed model; see][] {wooldridge2010}. The extension of our results to settings with growing $T$ is more complicated and beyond the scope of this paper, as this would require additional assumptions on the nature of the dependence across time.
\end{remark}

Let $\hat{\beta} = (\hat{\beta}_1, \hat{\beta}_2, ..., \hat{\beta}_K)'$ be the OLS estimator of $\beta$, that is
$$\hat{\beta} = \left(\sum_{n=1}^N {\bf x}_n {\bf x}_n'\right)^{-1}\sum_{n=1}^N {\bf x}_n y_n~.$$
In this paper, we focus on the use of $\hat{\beta}$ as a means of forming a test-statistic to perform inference on $\beta$. To that end, we study the asymptotic distribution of the following root:
$$T_k = \frac{\hat{\beta}_k - \beta_k}{\sqrt{\widehat{V}_{kk}}}~,$$
where $\hat{\beta}_k$ is the $k$th component of the OLS estimator of $\beta$ and $\hat{V}_{kk}$ is the $kk$th entry of an appropriate estimator of its asymptotic variance. For a specific value $\beta_{0k}$ of $\beta_k$ we call the resulting statistic the \emph{dyadic-robust $t$-statistic}. As mentioned in the introduction, the estimator of $\hat{V}$ we consider here is a sandwich estimator known in the literature as the dyadic-robust variance estimator:
$$\hat{V} = (\sum_{n=1}^N{\bf x}_n{\bf x}_n')^{-1}\Big(\sum_{n=1}^N\sum_{m=1}^N{\bf 1}_{n,m}\hat{u}_n\hat{u}_{m}{\bf x}_n{\bf x}_{m}'\Big)(\sum_{n=1}^N{\bf x}_n{\bf x}_n')^{-1}~,$$
where $\hat{u}_n = y_n - \hat{\beta}'{\bf x}_n$ and ${\bf 1}_{n,m}$ is an indicator function that equals 1 when $\psi(n) \cap \psi(m) \ne \emptyset$. 

Our goal is to specify conditions under which $T_k$ is asymptotically standard normal. As explained in the introduction, previous work on the dyadic-robust variance estimator either implicitly assumes an asymptotic normality result for $T_k$ \citep{fafgub2007, cammil2014}, or does not explore or employ such a result \citep{aronsam2015}. The contribution of this paper is to provide general conditions under which $T_k$ is asymptotically normal, and develop a degrees of freedom correction guided by our results.

\subsection{A Key Condition for our Central Limit Theorem}
To study the asymptotic distribution of $T_k$ we will employ a central limit theorem for dependency graphs proved in \cite{janson1988}. A key condition in the theorem, which we denote as Condition \hyperref[eq:janson]{(2)} in the appendix, plays a central role in our analysis. To simplify the exposition of our results, we introduce Condition \ref{cond:JC} below, which is a modified, but equivalent, condition. Remark \ref{rem:condequiv} in the appendix establishes the equivalence of these two conditions.

For each unit $g$, let $M_g$ be the number of dyads containing $g$. Recall that $N$ is the total number of \emph{dyads} in the data. Define 
$$\mathcal{M}^H = \max_g M_g, \hspace{4mm} \mathcal{M}^L = \min_g M_g~.$$ Note that by definition $\mathcal{M}^H \le G - 1$ and that $\frac{\mathcal{M}^LG}{2} \le N \le \frac{\mathcal{M}^HG}{2}$. 

\begin{condition}\label{cond:JC}  Let $\mathcal{M}^H$ be as above. Let $\sigma_N^2 = Var(\sum_n {\bf x}_n u_n)$. Given some additional assumptions (see Theorem \ref{thm:janson}), a sufficient condition for Janson's Theorem to apply in our framework is that
$$ \frac{(\frac{N}{\mathcal{M}^H})^{1/\ell}\mathcal{M}^H}{\sigma_N} \rightarrow 0 \hspace{2mm} \text{as} \hspace{2mm} N \rightarrow \infty~, $$
for some integer $\ell \ge 3$.
\end{condition}

Intuitively, Janson shows that in our framework the expression above gives a bound on the higher-order cumulants of the sequence of random variables we will study, and that these cumulants vanishing is sufficient to establish asymptotic normality.  A central theme of our paper is the fundamental connection Condition \ref{cond:JC} creates between $\mathcal{M}^H$ and $\sigma_N^2$. This connection will be made more clear throughout the rest of Section 2.

\begin{remark}\label{rem:Jansuff} It is important to emphasize that Condition \ref{cond:JC} is simply a \emph{sufficient} condition for asymptotic normality in our setting. Throughout Sections 2 and 3, we use Condition \ref{cond:JC} to motivate the assumptions we impose to ultimately prove our main asymptotic normality result (Theorem \ref{thm:main}). In Section 4 we perform a simulation study that explores what can happen when Condition \ref{cond:JC} does not hold. In Remark \ref{rem:CLTcompare} we comment on how our results would change if we were to use alternative central limit theorems for dependency graphs.
\end{remark}

\subsection{Asymptotic Frameworks}
\begin{figure}
\includegraphics[scale=0.350]{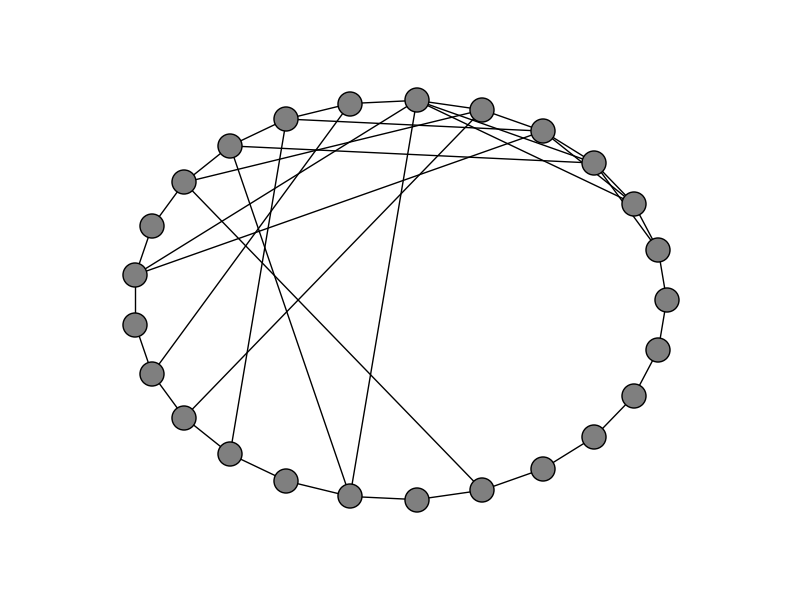} 
\includegraphics[scale=0.350]{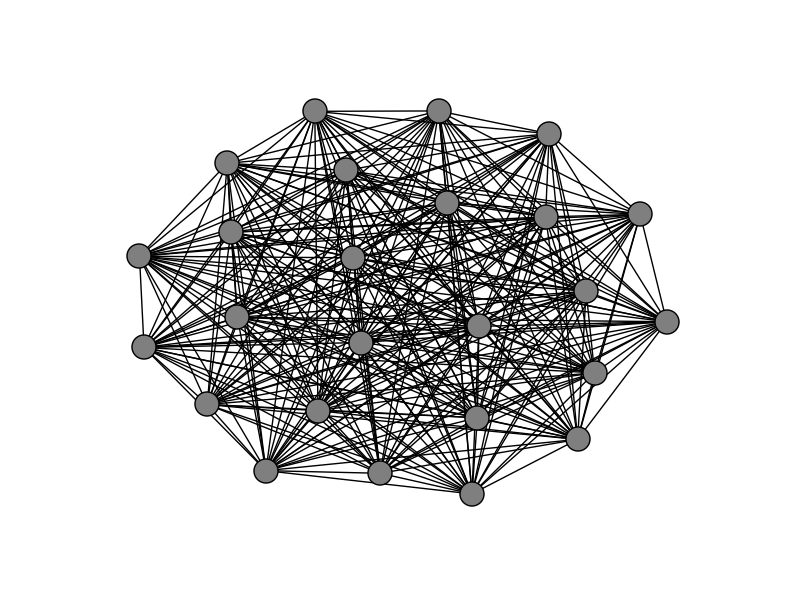}
\includegraphics[scale = 0.350]{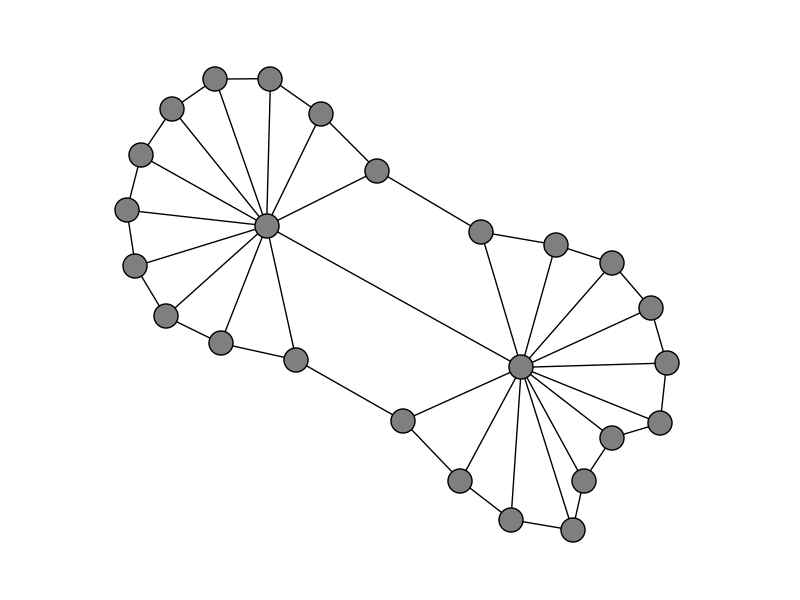}
\centering
\captionsetup{justification = centering} 
\caption{Clockwise from the top-left: Models S, D, and B with $G = 25$. \\
Units are the grey nodes, dyads are the black edges. \label{fig:diag}}
\end{figure}

We will now describe the two asymptotic frameworks that we consider in this paper. 
The first asymptotic framework we consider is one where $\mathcal{M}^H$ (and thus also $\mathcal{M}^L$) is bounded as $G \rightarrow \infty$.  This framework is relevant in settings where units have few links. Model S in Figure \ref{fig:diag} presents a configuration of the dyads in which we would expect such an approximation to be appropriate; note that there are 25 units (the grey nodes), but no units are contained in more than 6 dyads (the black edges). We will call this framework \hyperref[ass:AF1]{AF1}: 

\begin{assumption}\label{ass:AF1} ({AF1}) $\mathcal{M}^H < D$ for some constant $D$ as $G \rightarrow \infty$. \end{assumption} 
We will see in Section 3 that bounding $\mathcal{M}^H$ makes the analysis very clean, but \hyperref[ass:AF1]{AF1} may not be an appropriate framework for many settings of interest. In particular, \cite{cammil2014} suggest trade-data as an application where $\mathcal{M}^H$ could be very large relative to sample size. Model D in Figure \ref{fig:diag} presents a configuration of the dyads such that every unit is contained in a dyad with every other unit, which is the configuration employed in all of the simulation results presented in \cite{cammil2014} and \cite{aronsam2015}. Given such a configuration, we may not expect  \hyperref[ass:AF1]{AF1} to provide a good approximation in this setting for large $G$, hence we will also study an asymptotic framework where we allow $\mathcal{M}^H \rightarrow \infty$ as $G \rightarrow \infty$. This setting is more complicated, and several subtle issues will arise. 

The first issue that arises is that allowing $\mathcal{M}^H$ to grow as $G \rightarrow \infty$ now also frees $\mathcal{M}^L$ to grow as well. The most flexible possible framework would be one where we make no assumptions on $\mathcal{M}^L$ and simply require $\mathcal{M}^H \rightarrow \infty$. Unfortunately, in this case Janson's CLT cannot generally establish asymptotic normality when $\mathcal{M}^L$ is finite and fixed or grows too slowly. To illustrate, suppose that the error structure were of the form
$$u_{n(g,h)} = \alpha_g + \alpha_h + \epsilon_n~,$$
where $\alpha_g$, $\alpha_h$, $\epsilon_n$ are i.i.d. If $\mathcal{M}^L$ grows at rate $\log(G)$, and all but finitely many units form links at this rate, then $Var(\sum_{n}{\bf x}_nu_n)$ would grow at rate $N\log(G) \sim G(\log(G))^2$. If we suppose additionally that $\mathcal{M}^H$ grows at rate $G$, then Condition \ref{cond:JC} is \emph{not} satisfied, and Janson's CLT will not apply to our setting. The simulations of Section 4 will show that this has implications for inference when the number of dyads per unit varies wildly. In light of this, our second framework is given by:

\begin{assumption}\label{ass:AF2} ({AF2}) $\mathcal{M}^L \ge cG$ for some positive constant $c$. \end{assumption}
It is clear that this assumption is stronger than simply requiring that $\mathcal{M}^L \rightarrow \infty$, as it imposes an explicit rate of growth on $\mathcal{M}^L$ (and $\mathcal{M}^H$). It is possible to weaken this assumption slightly, but pinning down this rate clarifies our analysis and simplifies our notation. 
\begin{remark}\label{rem:conrem}
It is crucial to note that, although this framework allows every unit to be linked together, this does \emph{not} imply that every dyad can be correlated with every other dyad: consider the dense configuration of Model D presented in Figure \ref{fig:diag}. By construction, every dyad can be correlated with at most $46$ other dyads, despite there being $300$ dyads total. 
\end{remark}
Model B in Figure \ref{fig:diag} highlights the kind of configuration of the dyads alluded to above that this framework may have trouble capturing: note that most of the units are only contained in three dyads, but that two of the units are contained in many dyads. This type of configuration causes $\mathcal{M}^H$ to be very large relative to $\sigma_N$, which we have seen may cause Condition \ref{cond:JC} to fail. Configurations of this type can arise in empirical settings with dyadic data: as an example, the application in \cite{aronsam2015} features a configuration in which most of the units are contained in approximately $10$ dyads, but a handful of units are contained in upwards of $140$ dyads. Model B will play an important role in the simulation study of Section 4, as well as in developing our proposed degrees of freedom adjustment.

\clearpage

\subsection{The Rate of Growth of $Var(\sum_n{\bf x}_n u_n)$}
We have already seen that the connection between the asymptotic framework and the rate of growth of $Var(\sum_{n}{\bf x}_nu_n)$ plays an important role in establishing Condition \ref{cond:JC}. In this section we will highlight another important observation concerning the rate of growth of the variance that is essential to the results presented in Section 3.

First consider \hyperref[ass:AF1]{AF1}. We see immediately that imposing \hyperref[ass:AF1]{AF1} results in 
$$Var(\sum_{n=1}^N{\bf x}_nu_n) =  \sum_{n = 1}^N\sum_{m=1}^NCov({\bf x}_nu_n,{\bf x}_{m}u_{m})$$
growing at rate $N$. Hence under \hyperref[ass:AF1]{AF1} we will make the following assumption on the rate of growth:
\begin{assumption}\label{ass:rate1} $$\Omega := \lim_{G\rightarrow \infty} \frac{1}{N} Var(\sum_{n=1}^N{\bf x}_nu_n) = \lim_{G\rightarrow \infty}  \frac{1}{N} \sum_{n = 1}^N\sum_{m=1}^NCov({\bf x}_nu_n,{\bf x}_{m}u_{m}) \hspace{2mm} \text{is positive definite}~.$$
\end{assumption}

We have seen in Section 2.3 that under \hyperref[ass:AF2]{AF2}, when we allow $\mathcal{M}^H$ and $\mathcal{M}^L$ to grow as $G \rightarrow \infty$, the analysis can become more complicated. To clarify the exposition, we first consider a straightforward assumption that could be imposed under \hyperref[ass:AF2]{AF2}:
\begin{assumption}\label{ass:rate2a} $$\Omega := \lim_{G\rightarrow \infty} \frac{1}{NG}Var(\sum_{n=1}^N{\bf x}_nu_n) = \lim_{G\rightarrow \infty}  \frac{1}{NG} \sum_{n = 1}^N\sum_{m=1}^NCov({\bf x}_nu_n,{\bf x}_{m}u_{m}) \hspace{2mm} \text{is positive definite}~.$$
\end{assumption}
This assumption is completely analogous to the standard assumption made on the variance in clustered data with strong dependence, and is implicitly the assumption made in all of the results in \cite{aronsam2015}. Note that this assumption holds under an additive common shocks error specification:
$$u_{n(g,h)} = \alpha_g + \alpha_h + \epsilon_n~,$$
where $\alpha_g$, $\alpha_h$, $\epsilon_n$ are i.i.d, which is a specification considered in the simulations of both \cite{cammil2014} and \cite{aronsam2015}. More generally, Assumption \ref{ass:rate2a} would be appropriate in applications where the dependence results in a positive correlation between dyads.

It is important to point out, however, that \hyperref[ass:AF2]{AF2} does not \emph{imply} the rate of growth for $Var(\sum_n {\bf x}_n u_n)$ given in Assumption \ref{ass:rate2a}. For example, if instead the data were simply $i.i.d$, then $Var(\sum_{n}{\bf x}_nu_n)$ would grow at rate $N$ since $Cov({\bf x}_nu_n,{\bf x}_{m}u_{m}) = 0$ for $n \ne m$. In Section 4 we also consider a more interesting example: suppose we divide the units in the data into two groups, which we'll call $G_A$ and $G_B$, and specify the error term as
\[ u_{n(g,h)} = \begin{cases} 
      -(\alpha_g + \alpha_h) + \epsilon_n & \textrm{ if $g$ and $h$ belong to different groups} \\
      \alpha_g + \alpha_h + \epsilon_n & \textrm{ if $g$ and $h$ belong to the same group} \\
   \end{cases} \]
then by controlling the relative sizes of $G_A$ and $G_B$, we can achieve growth rates of the form $NG^r$ for any $r \in [0,1]$ in Model D while still maintaining the maximal amount of dependence (see the appendix for details). Such an error structure is a stylized example of a situation where a shock to a unit could affect certain dyadic relations positively, while affecting others negatively. These examples highlight the fact that Assumptions \ref{ass:dep} and \hyperref[ass:AF2]{AF2} can accommodate many plausible rates of growth of the variance. Since the goal of inference using robust-variance estimators is to be as agnostic as possible about the specific dependence structure in the data, we will consider these possibilities in our analysis. Our second assumption on the rate of growth of the variance under \hyperref[ass:AF2]{AF2} is given by:

\begin{assumption}\label{ass:rate2b} $$\Omega:=\lim_{G\rightarrow \infty} \frac{1}{NG^{r}} Var(\sum_{n=1}^N{\bf x}_nu_n) = \lim_{G\rightarrow \infty}  \frac{1}{NG^{r}} \sum_{n = 1}^N\sum_{m=1}^NCov({\bf x}_nu_n,{\bf x}_{m}u_{m}) \hspace{2mm} \text{is positive definite}~,$$
for some $r \in [0,1]$. \end{assumption} 
Note that Assumption \ref{ass:rate2a} is a special case of Assumption \ref{ass:rate2b} with $r = 1$.

%% file: results.tex
Recall that our goal is to study the asymptotic distribution of 
$$T_k = \frac{\hat{\beta}_k - \beta_k}{\sqrt{\widehat{V}_{kk}}}~.$$
Although \cite{fafgub2007} motivate the construction of $T_k$ by citing the work of \cite{conley1999}, asymptotic normality under our maintained assumptions does not follow from his results. Theorem \ref{thm:main} contains our main result. Our development proceeds in two steps: our first set of results establish that
$$\tau_N(\hat{\beta} - \beta)  \xrightarrow[]{d} N(0,V)~,$$
for some $V$ to be defined later in the section. It will turn out that the rate $\tau_N$ depends on both the specific characteristics of the dependence and the asymptotic framework considered. Our second set of results establish under which additional conditions we get that
$$\tau_N^2 \hat{V} \xrightarrow[]{p} V~.$$
By combining these two sets of results, we will see that $T_k$ is asymptotically standard normal under a range of assumptions.

For all of the results that follow, we make the following support assumption: 

\begin{assumption}\label{ass:support} $\{({\bf x}_n,u_n)\}_{n=1}^N$ have uniformly bounded support for all $N$. \end{assumption}

\begin{remark}\label{rem:bdsupp}
We choose to present the results under a bounded support assumption for expositional simplicity, but this assumption can be weakened at the expense of adding extra conditions on the moments of $({\bf x},u)$, without altering the substantive conclusions of the paper.  See the remarks after the proofs of Proposition \ref{prop:dist2} and Proposition \ref{prop:var2b} (Remarks \ref{rem:feller} and \ref{rem:moms}, respectively) for further discussion.
\end{remark}

\subsection{Asymptotic Normality}
Let us now study the asymptotic distribution of $\tau_N(\hat{\beta} - \beta)$ under both frameworks. Expanding the expression:
$$\tau_N(\hat{\beta} - \beta) = (\frac{1}{N}\sum_{n=1}^N{\bf x}_n{\bf x}_n')^{-1}\frac{\tau_N}{N}\sum_{n=1}^N{\bf x}_nu_n~.$$
We show that the first term converges in probability to $E[{\bf x}_n{\bf x}_n']^{-1}$ by showing that the variance of each component of $(1/N)\sum_{n=1}^N{\bf x}_n{\bf x}_n'$ converges to zero. We show that the second term converges in distribution to a normal by Janson's Theorem. First we state the result under \hyperref[ass:AF1]{AF1}:

\begin{proposition}\label{prop:dist1} Under Assumptions \hyperref[ass:AF1]{AF1}, \ref{ass:dep}, \ref{ass:rate1}, and \ref{ass:support}, 
$$\sqrt{N}(\hat{\beta} - \beta) \xrightarrow[]{d} N(0,V)~,$$
with $V = E({\bf x}_n{\bf x}_n')^{-1}\Omega E({\bf x}_n{\bf x}_n')^{-1}$, and $\Omega$ as in Assumption \ref{ass:rate1}. \end{proposition}
Note that under \hyperref[ass:AF1]{AF1} the rate $\tau_N$ is the standard $\sqrt{N}$. Intuitively, this is because \hyperref[ass:AF1]{AF1} imposes strong restrictions on the amount of dependence in the data. 

Our asymptotic normality result under \hyperref[ass:AF2]{AF2} is:

\begin{proposition}\label{prop:dist2} Under Assumptions \hyperref[ass:AF2]{AF2}, \ref{ass:dep}, \ref{ass:rate2b}  with $r > 0$ and \ref{ass:support},
$$\tau_N(\hat{\beta} - \beta) \xrightarrow[]{d} N(0,V)~,$$ 
where $\tau_N = \sqrt{\frac{N}{G^{r}}}$, $V = E({\bf x}_n{\bf x}_n')^{-1}\Omega E({\bf x}_n{\bf x}_n')^{-1}$, and $\Omega$ is as in Assumption \ref{ass:rate2b}.
\end{proposition}

Note that the rate of convergence $\tau_N$ is now scaled by the growth-rate of $Var(\sum_n {\bf x}_n u_n)$; the smaller is $r$, the closer we get to the standard $\sqrt{N}$ rate. Our result explicitly excludes the case $r = 0$, this is because with $r = 0$ and all of our imposed assumptions, Condition \ref{cond:JC} does not hold. In any case, the most likely situation in which $r$ is exactly zero would be if the data were in fact i.i.d, and asymptotic normality then follows from many standard CLTs.

Although it may seem problematic at first that $r$ will be unknown in practice, we will see in the following subsection that it is not necessary to know $r$ precisely in order to do inference on $\beta_k$.

\begin{remark}\label{rem:CLTcompare}
To prove Propositions \ref{prop:dist1} and \ref{prop:dist2} we employ a CLT for dependency graphs developed in \cite{janson1988}. However, we could also consider using other results, such as the Berry-Esseen bounds for dependency graphs developed in, for example, \cite{baldi1989}, \cite{penrose2003}, and \cite{chen2004}. For our purposes it seems that Janson's CLT is as general (if not more general) than other available results: for example, to prove Proposition \ref{prop:dist2} by employing the results in \cite{baldi1989} or \cite{penrose2003} would require us to assume that $r > 2/3$,  and the results in \cite{chen2004} do not allow us to establish an asymptotic normality result under AF2 and our maintained assumptions.
\end{remark}

\subsection{Consistency of $\hat{V}$}
We now turn to our second set of results: under what additional assumptions will the dyadic-robust variance estimator $\hat{V}$ converge to the asymptotic variance $V$? Recall that 
$$\hat{V} = (\sum_{n=1}^N{\bf x}_n{\bf x}_n')^{-1}\hat{\Omega}(\sum_{n=1}^N{\bf x}_n{\bf x}_n')^{-1}~,$$
\noindent where 
$$\hat{\Omega} = \Big(\sum_{n=1}^N\sum_{n'=1}^N{\bf 1}_{n,n'}\hat{u}_n\hat{u}_{n'}{\bf x}_n{\bf x}_{n'}'\Big)~.$$
To prove that $\tau_N^2 \hat{V} \xrightarrow[]{p} V$ we consider the ``bread'', $(\sum_{n}{\bf x}_n{\bf x}_n')^{-1}$, and the ``meat'', $\hat{\Omega}$, of the estimator separately. The convergence of $(1/N)\sum_{n}{\bf x}_n{\bf x}_n'$ to $E({\bf x}_n{\bf x}_n')$ was proved when deriving the asymptotic distribution of $\tau_N(\hat{\beta} - \beta)$. To show that $(\tau_N/N)^2\hat{\Omega} \xrightarrow[]{p} \Omega$, we show convergence in mean-square. As before, the result under \hyperref[ass:AF1]{AF1} is relatively straightforward:

\begin{proposition}\label{prop:var1} Under Assumptions \hyperref[ass:AF1]{AF1}, \ref{ass:dep}, \ref{ass:rate1}, and  \ref{ass:support}, we have that
$$N\hat{V} \xrightarrow[]{p} V~.$$
\end{proposition}

\noindent The result under \hyperref[ass:AF2]{AF2} will again need more qualifications. First we present the result under Assumption \ref{ass:rate2a}:

\begin{proposition}\label{prop:var2a} Under Assumptions \hyperref[ass:AF2]{AF2}, \ref{ass:dep}, \ref{ass:rate2a}, and \ref{ass:support}, we have that 
$$(N/G) \hat{V} \xrightarrow[]{p} V~.$$
\end{proposition}

Note that the rate has now slowed in accordance with the fact that the number of dependencies is growing for each unit. A special case of Proposition \ref{prop:var2a} for $\mathcal{M}^L = \mathcal{M}^H = G- 1$ appears in \cite{aronsam2015}, whose proof generalizes to ours. Our next result generalizes Proposition \ref{prop:var2a} to the setting where instead we impose Assumption \ref{ass:rate2b}:
\begin{proposition}\label{prop:var2b} Under Assumptions \hyperref[ass:AF2]{AF2}, \ref{ass:dep}, \ref{ass:rate2b} with $r > 1/2$, and \ref{ass:support}, we have that 
$$\tau_N^2 \hat{V} \xrightarrow[]{p} V~,$$
where $\tau_N$ is as in Proposition \ref{prop:dist2}.
\end{proposition}

Note that we have restricted the rate of growth of $Var(\sum_n {\bf x}_nu_n)$ even more than in Proposition \ref{prop:dist2}. We will provide some intuition as to why we require $r > 1/2$ in this result. In our proof we show convergence in mean-square, so heuristically we want to show that $Var((\tau^2/N^2)\hat{\Omega}) \rightarrow {\bf 0}$. It is the case that the variance of $\hat{\Omega}$ depends not only on the covariances between observations, but on their dependence in general. Now, the larger is $r$, the slower the growth of $\tau_N$, and hence the faster the convergence of $\tau^2/N^2$ to zero. For small values of $r$, the convergence of $\tau^2/N^2$ is too slow and cannot combat the growth in dependencies present in the data. This means that we cannot formally establish the consistency of $\hat{V}$ when the growth of $Var(\sum_n x_nu_n)$ is slow but the dependencies in the data are strong. We will study the implications of this via simulation in Section 4. 

\begin{remark}\label{rem:variid} Note that Propositions \ref{prop:var2a} and \ref{prop:var2b} do not establish the consistency of $\hat{V}$ when the data are i.i.d. As an aside, we prove this result separately in Proposition \ref{prop:var3} below.
\end{remark}
 
\begin{proposition}\label{prop:var3} Under Assumptions \hyperref[ass:AF2]{AF2}, \ref{ass:support} and the assumption that $\{({\bf x}_n,u_n)\}$ are i.i.d, we have that 
$$N\hat{V} \xrightarrow[]{p} V~.$$
\end{proposition}

\begin{remark}\label{rem:psd}
\cite{cammil2014} make the observation that there is no guarantee that $\hat{V}$ is positive semi-definite in finite samples. To account for this possibility they suggest the following modification: Consider the unitary decomposition $\hat{V} = U\Xi U'$ where $\Xi = diag[\lambda_1,...,\lambda_k]$ is the diagonal matrix of eigenvalues of $\hat{V}$. \cite{cammil2014} suggest replacing $\hat{V}$ by 
$$\hat{V}^{+} = U\Xi^{+}U'~,$$
where $\Xi^{+} = diag[\lambda_1^+, ..., \lambda_k^+]$, $\lambda_j^+ := max(\lambda_j,0)$ for all $j$. This modification guarantees that $\hat{V}^{+}$ is positive semi-definite. For the purposes of doing inference it would be even better if our estimator were positive-definite to guarantee that it is invertible. \cite{pol2011} suggests such a modification for a HAC estimator in a time-series setting, where $\widetilde{\lambda}_j^{+} := max(\lambda_j,\epsilon_n)$ with $\epsilon_n \rightarrow 0$ at a suitable rate, and proves its consistency there. Additionally, \cite{cammil2014} suggest a simple finite sample correction for $\hat{V}^+$ which may help alleviate some of its finite-sample bias. 
\end{remark}

\subsection{Asymptotic Normality of $T_k$}
With all of these results in hand, we can now state the main result of the paper:

\begin{theorem}\label{thm:main} Under Assumptions \ref{ass:dep}, \ref{ass:support}, and either:
\begin{itemize}
\item \hyperref[ass:AF1]{AF1} and Assumption \ref{ass:rate1},
\item \hyperref[ass:AF2]{AF2} and Assumption \ref{ass:rate2a},
\item \hyperref[ass:AF2]{AF2} and Assumption \ref{ass:rate2b} with $r > \frac{1}{2}$,
\end{itemize}
we have that,
$$T_k = \frac{\hat{\beta}_k  - \beta_k}{\sqrt{\widehat{V}_{kk}}} \xrightarrow[]{d} N(0,1)~.$$
\end{theorem}
\vspace{2mm}
\begin{remark}
With Theorem \ref{thm:main} in hand it is clear that a hypothesis test that uses the dyadic-robust $t$-statistic and the quantiles of a $N(0,1)$ distribution as critical values will be asymptotically valid. The same can be said for the coverage of an analogous confidence interval.
\end{remark}

We want to highlight two appealing aspects of this result. The first is that the limit distribution is the \emph{same} under both asymptotic frameworks. The second is that $T_k$ features a type of ``self-normalization" of $\tau_N$. Specifically, although the rate of convergence of the estimator depends on which asymptotic framework we consider and which assumptions we impose on the dependence, precise knowledge of this rate is not required to construct $T_k$. In this respect our paper is philosophically similar to \cite{hansen2007}, where he shows that the robust $t$-statistic in the linear panel model is asymptotically normal under various assumptions about the dependence across time. There is, however, one important caveat to our claim: under \hyperref[ass:AF2]{AF2} with Assumption \ref{ass:rate2b}, we only obtain the result when we impose that $r > 1/2$. This means that we have \emph{not} established an asymptotic normality result when the configuration of the dyads is dense and the dependence features many positively and negatively correlated dyads. We treat the results under Assumption \ref{ass:rate2b} as a form of robustness to small deviations from Assumption \ref{ass:rate2a}, which is the standard assumption imposed in the rest of the literature on strong dependence (i.e. clustering and its variations).  In Section 4, we will use the insights we have gained from our formal analysis of the problem to conduct a simulation study under a broader range of designs than those considered previously in the literature.

\begin{remark}\label{rem:comp_cam_aronow}
Theorem \ref{thm:main} is sufficiently general to accommodate the theoretical setup considered in \cite{aronsam2015}, as well as the simulation designs in  \cite{aronsam2015} and \cite{cammil2014}. In their settings, $\mathcal{M}_H = \mathcal{M}_L = G-1$, and $r = 1$, so that Theorem \ref{thm:main} applies under Assumptions \hyperref[ass:AF2]{AF2} and \ref{ass:rate2a}.
\end{remark}

\begin{remark}\label{rem:Mest}
\cite{cammil2014} and \cite{aronsam2015} also consider the use of the dyadic-robust variance estimator in settings with M-estimators or GMM estimators. Our results can be extended to these settings under a set of ``classical" regularity conditions \citep[see][Section 5.6]{van1998} which include, for example, the logit or probit models considered in \cite{cammil2014}. However, extending our results to M-estimation problems under more general regularity conditions  \citep[see for example][Theorem 5.23]{van1998} may require different tools, and hence different assumptions, than those considered in this paper. To this end, theoretical tools recently developed in \cite{lee2017} look promising.
\end{remark}

%% file: simulations.tex
In this section we perform a simulation-study to assess the accuracy of the normal approximation in finite samples, and the validity of the approximation for DGPs that do not satisfy the results presented in Section 3. In light of our results, we propose a novel degrees of freedom correction and study its performance via simulation as well. We use our formal results to guide the choice of designs. In an attempt to make the link between the simulations and our results as clear as possible, we consider simple designs.

The model we use throughout is a special case of what we have studied in Section 3:
$$y_n = 1 + \beta x_n + u_{n} ~,$$
with $x_n$ a scalar and $\beta = 0$. We consider two different specifications of $x_n$:
\begin{itemize}
\item When the $u_n$ are i.i.d, $x_n \sim U[0,1]$ i.i.d, so that $(x_n,u_n)$ are i.i.d.
\item When $u_n$ are not i.i.d, $x_{n(g,h)} = |z_g - z_h|$ where $z_g \sim U[0,1]$ i.i.d. for $g \in G$.
\end{itemize}

As pointed out in Remark \ref{rem:psd}, it is possible that $\hat{V}$ is not positive definite in finite samples, in particular when $G$ is small. For the simulations we perform, we use the estimator $\widetilde{V}$ which is constructed by an analogous construction to $\hat{V}^+$ presented in Remark \ref{rem:psd} but with $\tilde{\lambda}_j^+ = max(\lambda_j,\epsilon)$ for $\epsilon = 10^{-7}$. Whether or not we use $\widetilde{V}$, or $\hat{V}$ and drop those iterations where $\hat{V}$ is non-positive does not materially affect the results. In particular, for the sample sizes in which we get appropriate coverage, $\hat{V}$ and $\tilde{V}$ are always equal. Although we do not employ the finite sample corrections to $\hat{V}$ suggested in \cite{cammil2014}, we comment on them briefly in Section 4.3 below.

\subsection{Designs with Varying Levels of Sparsity}

First we will study the normal approximation under varying levels of ``sparsity" of the dyads relative to the number of units, while maintaining the standard variance Assumptions \ref{ass:rate1} and \ref{ass:rate2a}. We will consider three possible levels of sparsity:
\begin{itemize}
\item {\bf Model D (dense)} Every possible dyad is present in the data.
\item {\bf Model S (sparse)} A design where $\mathcal{M}^H = 5$ and $\mathcal{M}^L = 2$ regardless of sample size. 
\item {\bf Model B (both)} A design where $0.5 \log(G) \le \mathcal{M}^L \le 1.5 \log(G)$ and $0.5G \le \mathcal{M}^H \le G-1$, such that $Var(\sum_n x_n u_n)$ grows slower than required in Proposition \ref{prop:dist2}.
\end{itemize}

See the appendix for the specific construction of Models S and B. Figure \ref{fig:diag} in Section 2 was generated using these designs for $G = 25$. Given the results presented in Section 3, we expect Models D and S to behave well in large samples since they satisfy the conditions of our theorem. Model B does not satisfy the conditions of our theorem since the rate of growth of $\mathcal{M}^L$ is too slow. We saw in Section 2 that this has the potential to make $Var(\sum_n x_nu_n)$ grow too slowly for Janson's CLT to apply.

We will study these three models under two possible specifications of the error term. In the first specification, $u_n \sim U[-\sqrt{3},\sqrt{3}]$ i.i.d, so that $(x_n,u_n)$ are simply i.i.d. In the second specification we set:
$$u_{n(g,h)} = \alpha_g + \alpha_h + \epsilon_n~,$$ 
where $\alpha_g \sim U[-\sqrt{3},\sqrt{3}]$ i.i.d for $g = 1,2,...G$ and $\epsilon_n \sim U[-\sqrt{3},\sqrt{3}]$ i.i.d for $n = 1,2,...,N$. This is the additive common-shocks model discussed in Section 3. 

We perform 10,000 Monte Carlo iterations and record the number of times a $95\%$ confidence interval based on the normal approximation contains the true parameter $\beta = 0$. Table \hyperref[table:mainsim1]{1} presents the results under each specification. Recall that $G$ is the number of units in the data. 

\begin{table}[h]
\centering
\label{table:mainsim1}
\begin{tabular}{lllllll} \midrule \midrule
&&&&$G$&&\\
  &Specification                          & ~ 10 & ~ 25 & ~ 50 & ~100 & ~250 \\ \midrule
\multirow{3}{*}{Model D} & i.i.d &~69.8 & ~85.1   & ~91.3   & ~93.2  & ~94.4   \\
					&&(0.46)&(0.36)&(0.28)&(0.25)&(0.23)\\
                         & unit-level shock &  ~63.7  & ~86.2   & ~92.1   & ~93.6 & ~94.3    \\ 
                         &&(0.48)&(0.34)&(0.27)&(0.24)&(0.23)  \\
\multirow{3}{*}{Model S} & i.i.d &  ~70.0  & ~84.8   & ~90.1   & ~92.8    & ~94.5   \\
			&&(0.46) & (0.36) & (0.30) & (0.26) & (0.23) \\
                         & unit-level shock &  ~66.3  & ~82.8   & ~90.2   & ~92.9   & ~94.0    \\ 
                         && (0.47) &(0.38)&(0.30)&(0.26)&(0.24) \\
\multirow{3}{*}{Model B} & i.i.d &  ~69.6  & ~82.4   & ~86.9   & ~89.3    & ~93.1    \\
			&&(0.46)&(0.38)&(0.34)&(0.31)&(0.25) \\
                         & unit-level shock &  ~65.4 & ~81.0   & ~84.5   & ~86.1    & ~89.4   \\ 
                         && (0.48) & (0.39) & (0.36) & (0.35) &(0.31) \\ \bottomrule
\end{tabular}
\caption{Coverage percentages of a 95\% CI for $\beta = 0$, simulation SEs in parentheses.}
\end{table}

We note that for Models D and S, we get appropriate coverage in large samples. This is not surprising given our results, and is in line with the simulation results presented in \cite{cammil2014} and \cite{aronsam2015}. It is interesting to note that, despite their similar performances, there is significantly less data present in the simulations for a given $G$ under Model S than under Model D. For example, at $G = 50$ there are 1225 data points under Model D but only 86 data points under Model S. The fact that the approximations are comparable suggests that estimating the variance is much harder when the dyads are dense relative to the number of units. This makes sense seeing as the number of potential dependencies grows quadratically in Model D but is fixed in Model S. 

Given that it is the number of units that determines the accuracy of the approximation, it is also clear that we require many units to get adequate coverage. In many settings of interest (for example, social networks) the number of units required should not be insurmountable. In other settings such as international trade, the number of units required could be prohibitive. 

The final important observation we make is that the true coverage under Model B with unit-shocks, which was designed to fail the conditions of our theorem, is consistently worse than under either Models D and S, even when $G$ is relatively large. In fact, even with $G$ as large as $800$ (not formally reported) we do not see coverage higher than $92\%$ using this design. A similar phenomenon has been documented for the linear model with one-way clustering when cluster sizes differ ``wildly", in \cite{macweb2015} and \cite{cartersteig2013}. Because configurations in the spirit of Model B do arise in applications \citep[for example, the empirical application in][]{aronsam2015}, the degrees of freedom correction we propose in Section 4.3 is specifically designed to address this issue.

\subsection{Designs with Varying Growth Rates of $Var(\sum_n x_n u_n)$}

In this section we study the accuracy of a normal approximation under Assumption \ref{ass:rate2b}, which allows varying rates of growth of $Var(\sum_n x_n u_n)$. We noted in Section 3 that even though we cannot expect to know the exact rate in practice, using $T_k$ to test hypotheses does not require this knowledge. However, we also saw in Section 3 that when the model is dense, Proposition \ref{prop:var2b} does not establish consistency of $\hat{V}$ if the data is such that the growth rate of $Var(\sum_n x_nu_n)$ is too slow while still having many dependencies, which would be the case if the dependence features many positively and negatively correlated dyads.

Throughout this section we consider Models D and S. Our specification is constructed as follows: we divide the units in the data into two groups, which we call $G_A$ and $G_B$, and specify the error term as
\[ u_{n(g,h)} = \begin{cases} 
      -(\alpha_g + \alpha_h) + \epsilon_n & \textrm{ if $g$ and $h$ belong to different groups.} \\
      \alpha_g + \alpha_h + \epsilon_n & \textrm{ if $g$ and $h$ belong to the same group.} \\
   \end{cases} \]
Where $\alpha_g \sim U[-\sqrt{3},\sqrt{3}]$ i.i.d for $g = 1,2,...G$ and $\epsilon_n \sim U[-\sqrt{3},\sqrt{3}]$ i.i.d for $n = 1,2,...,N$.

By controlling the relative sizes of $G_A$ and $G_B$, we can achieve growth rates of the form $NG^r$ for any $r \in [0,1]$ in Model D while still maintaining the maximal amount of dependence (see the appendix for details). Although it is clear that this design is artificial, we think it is reasonable to consider situations where shocks at the unit-level can have differing effects across dyads.

As before, we perform 10,000 Monte Carlo iterations and record the number of times a $95\%$ confidence interval based on the normal approximation contains the true parameter $\beta = 0$. Table \hyperref[table:mainsim2]{2} presents the results for varying levels of $r$.

\noindent \begin{table}
\centering
\label{table:mainsim2}
\begin{tabular}{lllllll} \midrule \midrule
&$r$&&&$G$&&\\
     && ~ 10 & ~ 25 & ~ 50 & ~100 & ~250 \\ \midrule
\multirow{9}{*}{Model D} &$0$ & ~65.5  &  ~75.9  & ~80.4   &  ~83.6  &   ~84.0 \\ 
&&(0.48)&(0.43)&	(0.40) & (0.37) & (0.37) \\
 &$0.25$  & ~65.2    &  ~77.3  & ~82.8   &   ~86.7  &   ~89.8 \\
 &&(0.48)&(0.42)&(0.38)&(0.34)&(0.30) \\
 &$0.5$  & ~65.2    &  ~80.1  & ~87.3   &   ~91.5  &   ~94.7 \\
 &&(0.48)&(0.40)&(0.33)&(0.28)&(0.22)\\
& $0.75$  & ~65.3    &  ~82.6  & ~91.1   &   ~93.9  &   ~94.4 \\
&&(0.48)&(0.38)&(0.28)&(0.24)&(0.23) \\
& $1$  & ~63.7    &  ~82.6 & ~92.1   &  ~93.6  &  ~94.3 \\ 
&& (0.48) & (0.38) &  (0.27) & (0.24) & (0.23) \\ \midrule
\multirow{9}{*}{Model S} &$0$ & ~68.2    &  ~84.1  & ~90.3   &   ~92.8  &   ~93.8 \\
&&(0.47) & (0.37) & (0.30) & (0.26) & (0.24) \\ 
 &$0.25$  & ~68.6    &  ~84.1  & ~90.1   &  ~92.9  &   ~94.2 \\
 &&(0.46) & (0.37) & (0.30) & (0.26) & (0.23) \\
 &$0.5$  & ~68.6    &  ~84.0  & ~90.1   &  ~92.5  &  ~94.3 \\
 &&(0.46) & (0.37) & (0.30) & (0.26) & (0.23) \\
& $0.75$  & ~67.5    &  ~83.4  & ~90.4   &  ~92.6 &  ~94.1 \\
&&(0.47) & (0.37) & (0.29) & (0.26) & (0.24) \\
& $1$  & ~66.3    &  ~82.8  & ~90.2   &  ~92.9  &  ~94.0 \\
&& (0.47) & (0.38) & (0.30) & (0.26) & (0.24) \\ \bottomrule
\end{tabular}
\caption{Coverage percentages of a 95\% CI for $\beta = 0$. Simulation SEs in parentheses.}
\end{table}

Unsurprisingly given our results, coverage probabilities under Model S are at the nominal level in large samples. The results under Model D are more interesting: recall that Proposition \ref{prop:var2b} established consistency only for $r > 0.5$, and indeed we see that we get appropriate coverage in large samples for $r = 0.5, 0.75,$  and $r = 1$. In contrast, our simulations display poor coverage in large samples for $r = 0$ and $r = 0.25$, where Proposition \ref{prop:var2b} does not establish consistency. This suggests that $\hat{V}$ may be a poor approximation of the true asymptotic variance when there are a roughly equal number of negatively and positively correlated observations in the data. For this reason we feel that it is more appropriate to treat the results under Assumption \ref{ass:rate2b} as a form of robustness to small deviations from Assumption \ref{ass:rate2a}, which is the standard assumption imposed in the literature on strong dependence.

\subsection{A Degrees of Freedom Correction}

Given the simulation results in Section 4.1, and the fact that configurations like Model B do arise in empirical applications, we propose a new degrees of freedom correction to help guard against the potential for under-coverage. Instead of using the critical values from a $N(0,1)$ distribution to perform inference, we propose using the critical values from a $t_\kappa$ distribution where $\kappa$ is given by
$$\kappa = G\cdot\left(\frac{\text{med}_g\{M_g\}}{\mathcal{M}^H}\right)~,$$
where $\text{med}_g\{M_g\}$ denotes the median of the $\{M_g\}_{g=1}^G$. This degrees of freedom adjustment is similar in spirit to those proposed for analogous inference procedures using robust variance estimators in other settings: for example, \cite{belmac2002}, \cite{donald2007}, and \cite{imbkol2012} propose degrees of freedom adjustments for the heteroskedastic and clustered data settings \citep[for a textbook treatment, see][]{angrist2008}. Although it is ad-hoc, the intuition behind our choice of $\kappa$ is simple: when the configuration of the dyads satisfies our asymptotic normality assumptions, as in Models S and D, then $\text{med}_g\{M_g\}/\mathcal{M}^H$ is large, and hence the critical values derived from a $t_\kappa$ distribution approach the critical values derived from a $N(0,1)$ distribution for large $G$. On the other hand, for configurations like Model B where most of the units are contained in a few dyads but some units are contained in many dyads, we get that $\text{med}_g\{M_g\}/\mathcal{M}^H$ is very small, which results in a down-weighting of $\kappa$ and hence an enlargement of the critical value. In Table \hyperref[table:mainsim3]{3} we repeat the simulation exercise of Section 4.1, but with our degrees of freedom adjustment.  

Although we still see under-coverage for very unbalanced configurations such as Model B, we do a modest job of improving coverage in this setting while maintaining proper coverage in Models S and D; hence we see that a major benefit of employing this degrees of freedom correction is that it does not require the researcher to take a stand on whether or not their configuration is ``well-behaved". Moreover, this degrees of freedom correction could easily be implemented in any software package that computes dyadic standard errors and confidence intervals.  Our simulations suggest that, by combining our degrees of freedom correction with the finite sample corrections for $\hat{V}$ presented in \cite{cammil2014}, inference for most configurations of at least 150 units should behave as expected.

\begin{table}[h]
\centering
\label{table:mainsim3}
\begin{tabular}{lllllll} \midrule \midrule
&&&&$G$&&\\
  &Specification                          & ~ 10 & ~ 25 & ~ 50 & ~100 & ~250 \\ \midrule
\multirow{3}{*}{Model D} & i.i.d &~73.1 & ~86.4   & ~91.8   & ~93.5  & ~94.5   \\
					&&(0.44)&(0.34)&(0.27)&(0.25)&(0.23)\\
                         & unit-level shock &  ~67.6  & ~87.7   & ~92.8   & ~93.9 & ~94.5    \\ 
                         &&(0.47)&(0.33)&(0.26)&(0.24)&(0.23)  \\
\multirow{3}{*}{Model S} & i.i.d &  ~76.7  & ~87.7   & ~91.8   & ~93.4    & ~94.7   \\
			&&(0.43) & (0.33) & (0.27) & (0.25) & (0.22) \\
                         & unit-level shock &  ~73.2  & ~86.2   & ~91.6   & ~93.6   & ~94.2    \\ 
                         && (0.44) &(0.34)&(0.27)&(0.25)&(0.23) \\
\multirow{3}{*}{Model B} & i.i.d &  ~74.5  & ~88.3   & ~92.6   & ~94.5    & ~95.8    \\
			&&(0.46)&(0.32)&(0.26)&(0.23)&(0.2) \\
                         & unit-level shock &  ~69.5 & ~87.2   & ~91.7   & ~93.6    & ~93.4   \\ 
                         && (0.46) & (0.33) & (0.27) & (0.25) &(0.25) \\ \bottomrule
\end{tabular}
\caption{Coverage percentages of a 95\% CI for $\beta = 0$, with $t_\kappa$ critical vaues. Simulation SEs in parentheses.}
\end{table}

%% file: conclusion.tex
In this paper we have established a range of conditions under which the dyadic-robust $t$-statistic is asymptotically normal. We have also seen that in situations where our theorem does not apply, using a normal approximation of $T_k$ for inference may not be appropriate even for reasonable sample sizes. Our analysis suggests that, when we combine our degrees of freedom correction with the finite-sample corrections to $\hat{V}$ given in \cite{cammil2014}, inference should not be problematic for most datasets with roughly 150 units. However, if the data features a roughly equal number of positively and negatively correlated dyads, the dyadic-robust variance estimator may not provide a suitable approximation to the asymptotic variance. In our simulations this translated into confidence intervals that covered less often than the nominal coverage dictated, even in large samples. From a theoretical perspective it would be ideal to have a method of inference that is robust to this issue, but we expect that practitioners would consider such situations pathological in most settings of interest. 

In our opinion the most pressing issue to explore is that, as pointed out in Section 4, the normal approximation  of $T_k$ can be very poor when we do not have a lot of units in the data. A similar problem arises in the clustered-data setting, and recent papers have studied solutions for inference with few or even finitely many clusters \citep[see][for recent work in this area] {bch2011, cgm2008, crs2014, ibragimov2010, ibragimov2013}. A very promising option has been recently proposed in \cite{menzel2017}, who develops a bootstrap procedure for multiway/dyadic clustering that provides refinements whenever the limiting distribution is Gaussian. It would be interesting to see what other techniques could be adapted to the dyadic data setting as well.

\pagebreak

%% file: appendix.tex
\subsection{Proofs of the Results in Section 3}

In this section we will prove the results presented in section 3. To simplify the exposition, we first present the proofs in the special case of
$$y_n = \beta x_n + u_{n} ~,$$
where $x_n$ is a scalar. We then explain how to extend the proofs to the case where ${\bf x}_n$ is a vector (see Remarks \ref{rem:proofvec1} and \ref{rem:proofvec2}). 

A comment about the general strategy: to prove convergence in probability, we prove convergence in mean-square by finding an appropriate bound on the variance that converges to zero. To prove convergence in distribution to a normal, we use the following central limit theorem for dependency graphs proved in \cite{janson1988}. First we give the definition of a dependency graph for a family of random variables:

\begin{mydef} A graph $\Gamma$ is a dependency graph for a family of random variables if the following two conditions are satisfied:
\begin{itemize}
\item There exists a one-to-one correspondence between the random variables and the vertices of the graph.
\item If $\mathcal{V}_1$, $\mathcal{V}_2$ are two disjoint sets of vertices of $\Gamma$ such that no edge of $\Gamma$ had one endpoint in $\mathcal{V}_1$ and the other in $\mathcal{V}_2$, then the corresponding sets of random variables are independent. 
\end{itemize}
\end{mydef}
\noindent Now we state the theorem:

\begin{theorem}\label{thm:janson} {\bf (Janson 1988 Theorem 2)} Suppose that, for each $N$, $\{X_{Ni}\}_{i=1}^{N}$ is a family of bounded random variables; $|X_{Ni}| \le A_N$ a.s. Suppose further that $\Gamma_N$ is a dependency graph for this family and let $D_N$ be the maximal degree of $\Gamma_N$ (unless $\Gamma_N$ has no edges at all, in which case we set $D_N = 1$).
 Let $S_N = \sum_i X_{Ni}$ and $\sigma_N^2 = Var(S_N)$. If there exists and integer $\ell \ge 3$ such that
\begin{align} L_N:= \frac{(\frac{N}{D_N})^{1/\ell}D_NA_N}{\sigma_N} \rightarrow 0 \hspace{2mm} \text{as} \hspace{2mm} N \rightarrow \infty~, \label{eq:janson}\end{align}
then
$$\frac{S_N - E[S_N]}{\sigma_n} \xrightarrow[]{d} N(0,1) \hspace{2mm} \text{as} \hspace{2mm} N \rightarrow \infty~.$$
\end{theorem} 

\begin{remark}\label{rem:condequiv}As will be made clear in the proof of Proposition \ref{prop:dist2}, we use Theorem \ref{thm:janson} by definining an appropriate dependency graph for which $D_N = 2(\mathcal{M}^H-1)$, which establishes the equivalence between Condition \ref{cond:JC} and Equation \ref{eq:janson} for our purposes. Also note that although Janson's theorem applies to an \emph{array} of random variables, in the sense that for a given $i$, $X_{Ni}$ is allowed to change as $N$ grows, we do not use this feature for our results.
\end{remark}

Because the proofs under \hyperref[ass:AF1]{AF1} amount to special cases of the proofs under \hyperref[ass:AF2]{AF2}, we prove Propositions \ref{prop:dist2} and \ref{prop:var2b} before proving Propositions \ref{prop:dist1} and \ref{prop:var1} and \ref{prop:var2a}. 

\noindent {\bf Proof of Proposition \ref{prop:dist2}}
\begin{proof}
We have 
$$\tau_N(\hat{\beta} - \beta) = (\frac{1}{N}\sum_{n=1}^Nx_n^2)^{-1}\frac{\tau_N}{N}\sum_{n=1}^Nx_nu_n~.$$
First let's prove $((1/N)\sum_nx_n^2)^{-1} \xrightarrow[]{p} E[x_n^2]^{-1}$. By the continuous mapping theorem it is enough to show that $(1/N)\sum_nx_n^2 \xrightarrow[]{p} E[x^2]$. The expectation of $(1/N)\sum_n x_n^2$ is $E[x_n^2]$, so it suffices to show that 
$$Var(\frac{1}{N}\sum_{n=1}^Nx_n^2) \rightarrow 0~.$$
Expanding the variance gives 
$$Var(\frac{1}{N}\sum_{n=1}^Nx_n^2) = \frac{1}{N^2}\sum_{n=1}^N\sum_{m=1}^NCov(x_n^2,x_{m}^2)~.$$
For a fixed $n = n(g,h)$, we have at most $2\mathcal{M}^H - 1$ terms in the inner sum such that the covariance is nonzero. By Assumption \ref{ass:support}, $Cov(x_n^2,x_{m}^2)$ is uniformly bounded. Hence the sum over $n$ is of order $O(N\mathcal{M}^H)$. Thus under \hyperref[ass:AF2]{AF2} we have that $$Var(\frac{1}{N}\sum_{n=1}^Nx_n^2) = O\big(\frac{\mathcal{M}^H}{N}\big) = O\big(\frac{1}{G}\big) = o(1)~.$$
Next we'll prove that 
$$\frac{\tau_N}{N}\sum_{n=1}^Nx_nu_n \xrightarrow[]{d} N(0,\Omega)~.$$
We apply Janson's theorem to the family of random variables $\{x_nu_n\}_{n=1}^N$. A dependency graph $\Gamma_N = (\mathcal{V},E)$ for this family is the graph with vertex-set $\mathcal{V} = \{x_nu_n\}_{n=1}^N$, and edge set 
$$E = \{\{x_nu_n,x_{m}u_{m}\} : x_nu_n, x_{m}u_{m} \in V \hspace{2mm} \text{and} \hspace{2mm} \psi(n) \cap \psi(m) \ne \emptyset\}~.$$

That $\Gamma_N$ is a dependency graph for $\{x_nu_n\}_{n=1}^N$ follows immediately from Assumption \ref{ass:dep}. The maximal degree $D_N$ of $\Gamma_N$ is $2(\mathcal{M}^H - 1)$ by definition and, by Assumption \ref{ass:support}, $|x_nu_n| < A$ for all $N$ for some finite constant $A$. It remains to check Condition \hyperref[eq:janson]{(2)} of Janson's theorem. Let $\Omega_N = Var(\sum_nx_nu_n)$, then:
$$L_N = \frac{(\frac{N}{2(\mathcal{M}^H-1)})^{1/\ell}2(\mathcal{M}^H - 1)A}{\sqrt{\Omega_N}} = \frac{(\frac{N}{2(\mathcal{M}^H-1)})^{1/\ell}2(\mathcal{M}^H - 1)A}{\sqrt{NG^{r}}}\cdot\big(\frac{1}{NG^{r}}\Omega_N\big)^{-1/2}~.$$
Call the first term in the product $R_1$ and the second term $R_2$. $R_2$ converges to $\Omega^{-1/2}$ by Assumption \ref{ass:rate2b}. To evaluate the limit of $R_1$, we re-write everything in terms of the rates dictated by \hyperref[ass:AF2]{AF2}, which gives us
$$R_1 = O\Big(\frac{G^{1/\ell}}{G^{r/2}}\Big)~.$$
Choose $\ell$ sufficiently large such that $\frac{1}{\ell} - \frac{r}{2} < 0$, which is possible since $r > 0$. It then follows that $R_1 \rightarrow 0$. Now that we have established Condition \hyperref[eq:janson]{(2)} of Janson's theorem, we have that
$$\frac{\sum_{n=1}^Nx_nu_n}{\sqrt{\Omega_N}} \xrightarrow[]{d} N(0,1)~.$$
Re-writing:
$$\frac{\sum_{n=1}^Nx_nu_n}{\sqrt{\Omega_N}} = \big(\frac{1}{\sqrt{NG^{r}}}\sum_{n=1}^Nx_nu_n\big)\cdot(\frac{1}{NG^{r}}\Omega_N)^{-1/2} = \frac{\tau_N}{N}\sum_{n=1}^Nx_nu_n\cdot(\frac{1}{NG^{r}}\Omega_N)^{-1/2}$$
It thus follows by Assumption \ref{ass:rate2b} that 
$$\frac{\tau_N}{N}\sum_{n=1}^Nx_nu_n \xrightarrow[]{d} N(0,\Omega)~.$$
Applying Slutsky's Theorem to
$$\tau_N(\hat{\beta} - \beta) = (\frac{1}{N}\sum_{n=1}^Nx_n^2)^{-1}\frac{\tau_N}{N}\sum_{n=1}^Nx_nu_n ~,$$
we can conclude that
$$\tau_N(\hat{\beta}-\beta) \xrightarrow[]{d} N(0,V) ~,$$
as desired.
\end{proof}

\begin{remark}\label{rem:feller} As noted after the statement of Assumption \ref{ass:support}, we could weaken the uniform boundedness assumption by using Janson's theorem with a Lindeberg-type condition (see Remark 3 in \cite{janson1988}). For example, if we instead assume that $x_nu_n$ has bounded $2 + \delta$ moments, then the result could be proved under the additional assumption that $r > 2/(2+\delta)$, which agrees with our result as we take $\delta$ to infinity. \end{remark}

\begin{remark}\label{rem:proofvec1}  The general case is proved as follows: to show the convergence of  $(1/N)\sum_n{\bf x}_n{\bf x}_n'$ we repeat the argument above but component-wise. To show the convergence of $(\tau_N/N)\sum_n{\bf x}_nu_n$ we use Janson's Theorem paired with the Cramer-Wold device. \end{remark}

\noindent {\bf Proof of Proposition \ref{prop:dist1}}
\begin{proof}
This proof follows by a similar argument to the proof of Proposition \ref{prop:dist2}. We will sketch it here. Expanding:
$$\sqrt{N}(\hat{\beta} - \beta) = (\frac{1}{N}\sum_{n=1}^Nx_n^2)^{-1}\frac{1}{\sqrt{N}}\sum_{n=1}^Nx_nu_n~.$$
The first term of the product converges to $E(x^2)^{-1}$ by the same argument as above. The second term of the product converges to a normal by an application of Janson's theorem \ref{thm:janson} where we note that now under \hyperref[ass:AF1]{AF1} the maximal degree $D_N = 2(\mathcal{M}^H-1)$ of the dependency graph is bounded. Hence by Janson's theorem we have that
$$\frac{\sum_{n=1}^Nx_nu_n}{\Omega_N} \xrightarrow[]{d} N(0,1)~.$$
By a similar calculation to the one done in the proof of Proposition \ref{prop:dist2}, we get that
$$\frac{1}{\sqrt{N}}\sum_{n=1}^Nx_nu_n \xrightarrow[]{d} N(0,\Omega)~,$$
and thus Proposition \ref{prop:dist1} follows by an application of Slutsky's theorem. 
\end{proof}

\noindent{\bf Proof of Proposition \ref{prop:var2b}}
\begin{proof}
We follow the general strategy of \cite{aronsam2015}. First we introduce some notation for the proof. Given a dyadic index $n$, define $\tilde{n} = \psi(n)$. Recall that
$$\hat{V} = (\sum_{n=1}^Nx_n^2)^{-1}\hat{\Omega}(\sum_{n=1}^Nx_n^2)^{-1}~,$$
\noindent where 
$$\hat{\Omega} = \sum_{n=1}^N\sum_{m=1}^N{\bf 1}_{n,m}\hat{u}_n\hat{u}_{m}x_nx_{m}~.$$
We proved in Proposition \ref{prop:dist2} that $\big((1/N)\sum_{n=1}^Nx_n^2\big)^{-1} \xrightarrow[]{p} E(x^2)^{-1}$. So it remains to show that
$$\frac{\tau_N^2}{N^2}\hat\Omega \xrightarrow[]{p} \Omega~.$$
Re-writing $\hat{u}_n = u_n - (\hat{\beta} - \beta)x_n$ and expanding gives:
$$\frac{\tau_N^2}{N^2}\hat\Omega = \frac{\tau_N^2}{N^2}\Big(R_1 + R_2 + R_3 + R_4\Big)~,$$
where
$$R_1 = \sum_{n=1}^N\sum_{m=1}^N{\bf 1}_{n,m}x_nx_{m}u_nu_{m}~,$$
$$R_2 = - \sum_{n=1}^N\sum_{m=1}^N{\bf 1}_{n,m}x_nx_{m}^2u_n(\hat{\beta} - \beta)~,$$
$$R_3 = -  \sum_{n=1}^N\sum_{m=1}^N{\bf 1}_{n,m}x_n^2x_{m}u_{m}(\hat{\beta} - \beta)~, \hspace{2mm} \text{and}$$
$$R_4 =  \sum_{n=1}^N\sum_{m=1}^N{\bf 1}_{n,m}x_n^2x_{m}^2(\hat{\beta} - \beta)^2~.$$
We will show that $(\tau_N^2/N^2)R_1 \xrightarrow[]{p} \Omega$ while the rest converge in probability to zero.
As usual, we show convergence in mean-square. For the first term, it is the case by definition that
$$\lim_{G\rightarrow\infty}E\Big[\frac{\tau_N^2}{N^2}R_1\Big] = \Omega~,$$
so it suffices to show that
$$\lim_{G\rightarrow\infty}Var\Big(\frac{\tau_N^2}{N^2}R_1\Big) = 0~.$$
Expanding:
$$Var\Big(\frac{\tau_N^2}{N^2}R_1\Big) = \frac{\tau_N^4}{N^4}\Big(\sum_{i}\sum_{j}\sum_{k}\sum_{l}Cov({\bf 1}_{i,j}x_{i}x_{j}u_{i}u_{j},{\bf 1}_{k,l}x_{k}x_{l}u_{k}u_{l})\Big)~.$$
By Assumption \ref{ass:support}, the summands are uniformly bounded, so in order to get a suitable bound on the sum we need to understand under what conditions  
$$Cov({\bf 1}_{i,j}x_{i}x_{j}u_{i}u_{j},{\bf 1}_{k,l}x_{k}x_{l}u_{k}u_{l}) \ne 0~.$$ First it is clear that we must have 
\begin{align} \tilde{i} \cap \tilde{j} \ne \emptyset \hspace{2mm} \text{and} \hspace{2mm} \tilde{k} \cap \tilde{l} \ne \emptyset~. \label{eq:vcond1}\end{align}
Given \hyperref[eq:vcond1]{(2)} holds, expand the covariance:
$$Cov(x_{i}x_{j}u_{i}u_{j},x_{k}x_{l}u_{k}u_{l}) = E[x_{i}x_{j}u_{i}u_{j}x_{k}x_{l}u_{k}u_{l}] - E[x_{i}x_{j}u_{i}u_{j}]E[x_{k}x_{l}u_{k}u_{l}]~,$$
then we see that we must also have that
\begin{align} \tilde{i} \cap \tilde{k} \ne \emptyset \hspace{2mm} \text{or} \hspace{2mm} \tilde{i} \cap \tilde{l} \ne \emptyset \hspace{2mm} \text{or} \hspace{2mm} \tilde{j} \cap \tilde{k} \ne \emptyset \hspace{2mm} \text{or} \hspace{2mm} \tilde{j} \cap \tilde{l} \ne \emptyset~. \label{eq:vcond2}\end{align}
Let $S$ be the set of tuples $(i,j,k,l) \in N^4$ such that conditions \hyperref[eq:vcond1]{(2)} and \hyperref[eq:vcond2]{(3)} hold, then the cardinality of $S$, denoted $|S|$, is an upper bound on the number of nonzero terms in the sum. Our goal is to find an upper-bound on $|S|$.

Fix an index $i$, then there are $O(\mathcal{M}^H)$ indices $j$ such that $\tilde{i} \cap \tilde{j} \ne \emptyset$. Now, for a fixed $i$ and $j$ such that $\tilde{i} \cap \tilde{j} \ne \emptyset$, there are $O(\mathcal{M}^H)$ possible indices $k$ such that $\tilde{i} \cap \tilde{k} \ne \emptyset$ or $\tilde{j} \cap \tilde{k} \ne \emptyset$. For a fixed $i$, $j$ and $k$ such that the above hold, there are $O(\mathcal{M}^H)$ possible indices $l$ such that $\tilde{k} \cap \tilde{l} \ne \emptyset$.

Similarly, for a fixed $i$ and $j$ such that $\tilde{i} \cap \tilde{j} \ne \emptyset$, there are $O(\mathcal{M}^H)$ possible indices $l$ such that $\tilde{i} \cap \tilde{l} \ne \emptyset$ or $\tilde{j} \cap \tilde{l} \ne \emptyset$. For a fixed $i$, $j$ and $l$ such that the above hold, there are $O(\mathcal{M}^H)$ possible indices $k$ such that $\tilde{k} \cap \tilde{l} \ne \emptyset$.

Thus there are $N\cdot O(\mathcal{M}^H) \cdot O(\mathcal{M}^H) \cdot O(\mathcal{M}^H) = O(N(\mathcal{M}^H)^3)$ possible indices $i,j,k,l$ such that $(i,j,k,l) \in S$. Re-writing using the rates dictated by \hyperref[ass:AF2]{AF2} gives us that $|S| = O(G^5)$ and that
$$\frac{\tau_N^4}{N^4} \le K\frac{1}{G^{4+2r}}~,$$ 
for some positive constant K. Therefore we can conclude that
$$Var\big(\frac{\tau_N^2}{N^2}R_1\big) \le \frac{1}{G^{4+2r}}O(G^5) = o(1)$$
for $r > \frac{1}{2}$.
Thus we have shown that $(\tau_N^2/N^2)R_1 \xrightarrow[]{p} \Omega$. Next, we must show that the remaining terms converge to $0$ in probability. All three terms follow by similar arguments so we will only present the argument for $R_2$. We wish to show that 
$$\frac{\tau_N^2}{N^2}\Big(\sum_{n=1}^N\sum_{m=1}^N{\bf 1}_{n,m}x_nx_{m}^2u_n(\hat{\beta} - \beta)\Big) \xrightarrow[]{p} 0~.$$
We know from Proposition \ref{prop:dist2} that $\tau_N^{1-\epsilon}(\hat{\beta} - \beta) = o_P(1)$ for any $\epsilon > 0$, so it suffices to show that 
$$\frac{\tau_N^{1+\epsilon}}{N^2}\Big(\sum_{n=1}^N\sum_{m=1}^N{\bf 1}_{n,m}x_nx_{m}^2u_n\Big) = O_P(1)~,$$
for some $\epsilon > 0$. Note that from Assumption \ref{ass:support} and the definition of $\tau$,
$$E\Big[\frac{\tau_N^{1+\epsilon}}{N^2}\Big(\sum_{n=1}^N\sum_{m=1}^N{\bf 1}_{n,m}x_nx_{m}^2u_n\Big)\Big] \rightarrow 0~,$$
for $\epsilon$ sufficiently small, and that 
$$Var\Big(\frac{\tau_N^{1+\epsilon}}{N^2}\Big(\sum_{n=1}^N\sum_{m=1}^N{\bf 1}_{n,m}x_nx_{m}^2u_n\Big)\Big) \rightarrow 0~,$$
which can be shown by similar arguments to what we have done above. Similarly, the third and fourth terms also converge to zero in probability, and hence we have that
$$\frac{\tau_N^2}{N^2}\hat{\Omega} \xrightarrow[]{p} \Omega~,$$
and ultimately that
$$\tau_N^2\hat{V} \xrightarrow[]{p} V$$
as desired.
\end{proof}
\begin{remark}\label{rem:moms} Note that this proof relied on the bounded support Assumption \ref{ass:support} to ensure that the covariance terms in the summands were uniformly bounded. In general we would need to make assumptions on these covariances if we were to weaken the support assumption used here.\end{remark}
\begin{remark}\label{rem:proofvec2}  The general case is proved as follows: to show the convergence of $R_1$ to $\Omega$ we repeat the argument above but component-wise. To show the convergence of $R_2$, $R_3$ and $R_4$ to zero we modify the argument slightly. Consider
$$R_2 = \sum_{n=1}^N\sum_{m=1}^N{\bf 1}_{n,m}{\bf x}_n{\bf x}_{m}'(\hat{\beta} - \beta)'{\bf x}_nu_n~.$$
To show $R_2 \xrightarrow[]{p} {\bf 0}$ we will show that $||R_2|| \xrightarrow[]{p} 0$ where $||\cdot||$ is the Frobenius norm. Using the triangle inequality, the matrix Schwartz Inequality, and the definition of the Frobenius norm we get that 
$$||R_2|| \le \sum_{n=1}^N\sum_{m=1}^N{\bf 1}_{n,m}||{\bf x}_n||^2||{\bf x}_{m}||\cdot||\hat{\beta} - \beta||\cdot|u_n|~.$$
The result then follows from the arguments presented above.
\end{remark}

\noindent{\bf Proof of Proposition \ref{prop:var1}}
\begin{proof} Again we follow the same strategy as the proof of Proposition \ref{prop:var2b}. Now when getting a bound on the variance of $R_1$, we have that, for each fixed $i$, there are only finitely many $j,k,l$ such that $(i,j,k,l) \in S$. Therefore there are $O(N)$ nonzero terms in the sum. It then follows that the variance of $(1/N)R_1$ converges to zero and the proof goes through in the same manner as before.
\end{proof}
\noindent{\bf Proof of Proposition \ref{prop:var2a}}
\begin{proof} This is just a special case of Proposition \ref{prop:var2b}.
\end{proof}
\pagebreak
\noindent {\bf Proof of Proposition \ref{prop:var3}}
\begin{proof} We follow the same strategy as the proof of Proposition \ref{prop:var2b}. Now when getting a bound on the variance of $R_1$, we have that the terms in the sum are nonzero if and only if $i = j= k =l$, so that there are $N$ nonzero terms in the sum. It then follows that the variance of $(1/N)R_1$ converges to zero and the proof goes through in the same manner as before.
\end{proof}

\subsection{Simulation Design Details}
In this section we provide some details about the construction of our designs not mentioned in the main body. All of the simulations in the paper were performed using numpy in Python 2.7.

\noindent {\bf Construction of Model S: }Dyads were included in Model S according to the following rule:
\begin{itemize}
\item $(g,h)$ is included if $|g - h| = 1$
\item $(1,G)$ is included
\item $(g, 2g)$ is included if $g \le \left \lfloor{\frac{G}{2}}\right \rfloor$
\item $(g,3g)$ is included if $g \le \left \lfloor{\frac{G}{3}}\right \rfloor$
\end{itemize} 

\noindent {\bf Construction of Model B:} Dyads were included in Model B according to the following rule:
\begin{itemize}
\item $(g,h)$ is included if $|g - h| = 1$ and $g, h < G - 1$
\item $(1,G-2)$ is included
\item If $G = 100$, $G = 250$, or $G = 800$, $(g,h)$ is included if $|g-h| = 2$ and $g,h < G -1$
\item If $G = 100$ or $G =  250$, or $G = 800$, $(1,G-3)$ and $(2,G-2)$ are included
\item If $G = 800$, $(g,h)$ is included if $|g - h| \le 4$ and $g, h < G - 1$
\item If $G = 800$, $(1,G-4)$, $(1,G-5)$, $(2,G-3)$, and $(2,G-4)$ are included
\item $(g, G-1)$ is included for all $g \le \left \lfloor{\frac{G}{2}}\right \rfloor$, and $(g, G)$ is included for all $g > \left \lfloor{\frac{G}{2}}\right \rfloor$

\end{itemize}

\noindent {\bf Construction of $G_A$ and $G_B$}: Recall that in Section 4.2 we constructed our design by partitioning the units into two groups $G_A$ and $G_B$ and then specifying the error term as
\[ u_{n(g,h)} = \begin{cases} 
      -(\alpha_g + \alpha_h) + \epsilon_n & \textrm{ if $g$ and $h$ belong to different groups.} \\
      \alpha_g + \alpha_h + \epsilon_n & \textrm{ if $g$ and $h$ belong to the same group.} \\
   \end{cases} \]
Where $\alpha_g \sim U[-\sqrt{3},\sqrt{3}]$ i.i.d for $g = 1,2,...G$ and $\epsilon_n \sim U[-\sqrt{3},\sqrt{3}]$ i.i.d for $n = 1,2,...,N$.
In order to achieve a rate of growth of $NG^r$ for $Var(\sum_n {\bf x}_n u_n)$, we construct $G_A$ and $G_B$ as follows:
$$G_A = \Big\{g \in G : g \le  \left \lfloor{\frac{G - G^s}{2}}\right \rfloor\Big\}~,$$
$$G_B = \Big\{g \in G : g >  \left \lfloor{\frac{G - G^s}{2}}\right \rfloor\Big\}~,$$
where $s = \frac{1 + r}{2}$.
Expanding $Var(\sum_n {\bf x}_n u_n)$ shows that it indeed grows at rate $NG^r$.
\pagebreak